\documentclass[12pt]{amsart}
\usepackage{amscd,amsthm,amsmath} % SRC Specials: DVI [Inverse] Search
\newtheorem{thm}{Theorem}[section]
\newtheorem{cor}[thm]{Corollary}
\newtheorem{lem}[thm]{Lemma}
\newtheorem{prop}[thm]{Proposition}

\newtheorem{ques}[thm]{Question}

\newtheorem{claim}[thm]{Claim}
\theoremstyle{definition}
\newtheorem{defn}[thm]{Definition}
\theoremstyle{remark}
\newtheorem{rem}[thm]{Remark}
\numberwithin{equation}{section}
\begin{document}

\title{Having the H-space structure is not a generic property}%
\author{Jianzhong Pan }%
\address{Institute of Math.,Academia Sinica ,Beijing 100080, China }%
\email{pjz@math03.math.ac.cn}%
\thanks{The  author is partially supported by the NSFC project
19701032 and ZD9603 of Chinese Academy of Science }%

\date{June,5,2000}%
%\dedicatory{}%
%\commby{}%
% ----------------------------------------------------------------
\begin{abstract}
 In this note, we answer in negative  a question posed
 by McGibbon\cite{mcgib} about the generic property of H-space structure.
 In fact we verify  the conjecture of
 Roitberg \cite{roit}. Incidentally, the same example also answers
 in negative the  open problem 10 in McGibbon\cite{mcgib}
\end{abstract}
\maketitle
% ----------------------------------------------------------------
 \section{Introduction} \label{S:intro}
Let $X$ be a connected CW complex, the L-S category of $X$,
$cat(X)$, of $X$ is the least integer $k\geq 0$ such that $X$ can
be covered by $k+1$ open subsets which are contractible in $X$. Of
course the condition that $X$ is a CW complex is unnecessary for
the above definition to have a meaning . However it is in this
context that a rich theory of category exists. Recent works in
rational homotopy theory gave rise to a theory of rational L-S
category. This makes it possible to calculate the rational L-S
category and thus attack the rational Ganea
Conjecture\cite{ganea}. On the other hand works by N.Iwase,
H.Scheerer, and D.Stanley provided the method to determine the L-S
category itself in some case which lead to the construction of
counterexamples to the Ganea Conjecture, see, e.g., \cite{stan}.
Besides the application to the Ganea Conjecture , another
interesting application of these ideas was given by
Roitberg\cite{roit}. To explain this we have to state one problem
posed by McGibbon\cite{mcgib}:
\begin{ques}\label{T:ques1}
If $X$ and $Y$ have the same Mislin genus ,i.e. $X_{(p)}\simeq
Y_{(p)}$ for all $p$ where $X_{(p)}$ is the $p$-localization of
$X$, does it follow that $cat(X)=cat(Y)$?
\end{ques}
In his paper\cite{roit}, applying some results about category by
Iwase\cite{iwase} and results about phantom maps , Roitberg was
able to answer the above question negatively. His main result can
be stated as follows
\begin{thm}
Let $\phi:\Sigma K(\mathbb{Z},5)\to S^4$ be an essential, special,
phantom map and $X$ be the mapping cone of $\phi$ .Then
$cat(X)=2$.
\end{thm}
\begin{rem}
It is well known that $cat(X)=1$ iff $X$ is a co-H-space . It
follows from the above theorem that $X$ is not a co-H-space . On
the other hand it easy to know that $S^4 \vee \Sigma^2 K(Z,5)$ has
the same Mislin genus with $X$ and is a co-H-space .
\end{rem}
From this point of view , the theorem above answers in negative
the following
\begin{ques}\label{T:ques2}
If  $X$ and $Y$ have the same Mislin genus and $X$ is co-H-space ,
does it follow that $Y$ is also a co-H-space?
\end{ques}
Question\ref{T:ques1} has an obvious Eckmann-Hilton dual . However
it may not be a good question at present time since the dual L-S
category is not well developed. A manageable problem is the
obvious dual of Question\ref{T:ques2} which has  been posed by
McGibbon in \cite{mcgib}. A more precise conjecture  was given by
Roitberg in his paper\cite{roit}. The purpose of this paper is to
establish Roitberg's conjecture and thus also answers in negative
McGibbon's problem.
\begin{thm}\label{T:main}
Let $\psi:K(\mathbb{Z},2) \to \Omega S^{6}$ be an essential,
special, phantom map and $Z$ be the homotopy fiber of $\psi$. Then
$Z$ and $K(\mathbb{Z},2)\times \Omega^2 S^{6}$ have the same
Mislin genus and $Z$ is not an H-space.
\end{thm}
Since the dual L-S category is not well developed and Iwase's
paper\cite{iwase} is unavailable to the author,  method from the
well developed theory of H-space will be used in stead. Another
feature is the application of the newly developed Gray index of
phantom map \cite{lestrom},\cite{mstrom}. Actually , by duality ,
our method also gives an alternative proof of the main results in
\cite{roit}.

Incidentally , the example constructed above combined with Theorem
3.4 in \cite{roitharp} also provides a negative answer to the open
problem 10 in \cite{mcgib}: Is $X$  an H-space if each of its
Postnikov approximations $X^{(n)}$ is ?  Actually we have the
following
\begin{thm}
Let  $\psi:K(\mathbb{Z},2) \to \Omega S^{6}$ be an essential,
special, phantom map and $Z$ be the homotopy fiber of $\psi$. Then
$Z^{(n)}$ is an H-space for each n but $Z$ is not.
\end{thm}
\begin{proof}
It follows immediately from Theorem 3.4 in \cite{roitharp} that
$Z^{(n)}$ and $(K(\mathbb{Z},2)\times \Omega^2 S^{6})^{(n)}$ have
the same homotopy type and thus $Z^{(n)}$ is  an H-space for each
n.
\end{proof}

In this paper all spaces involved are assumed to be 1-connected CW
complexes with finite type.

The author would like to thank Prof. Roitberg for his interest in
this work and for pointing out a fatal error in the earlier
version of this paper . It is to correct that error that we find
the application of Gray index to this work.
\section{Background about H-spaces and Phantom maps}
First we will recall some backgrouds about H-spaces , see
\cite{az} for details. An H-space is a space $X$ with a map $\mu:X
\times X \to X$ such that $\mu|_{X\vee X} = F$ where $F:X\bigvee X
\to X$ is the natural folding map. An H-map between H-spaces is a
map of spaces $f:X \to Y$ such that the following diagram commutes
up to homotopy.
\[
\begin{CD}
X \times X @>f\times f>> Y \times Y \\
 @V\mu_X VV  @V\mu_Y VV \\
 X  @>f>> Y
\end{CD}
\]
In this case we say that $f$ is a $\mu_X - \mu_Y $ H-map.
 Two elementary but important results are the followings
\begin{prop}
Let $(X,\mu_X)$ be an H-space. Then, for any space M, $[M,X]$ is
an algebraic loop, i.e., for any $f,g \in [X,Y]$ there exists a
unique $D_{f,g} \in [M,X]$ such that
\[
\mu_*(D_{f,g} , g) = f
\]
\end{prop}
\begin{prop}\label{T:fiber}
If $f:(X,\mu_X)\to (Y,\mu_Y)$ is an H-map , then the homotopy
fiber of $f$ is an H-space.
\end{prop}
Thus it is important to know when is a map an H-map or what is the
obstruction for a map to be an H-map.
\begin{defn}
Let $(X,\mu)$ and  $(Y,\mu')$ be H-spaces and $f: X \to Y$ be a
map of spaces. H-derivation of $f$ is the map
\[
HD(f)\in [X\wedge X , Y]
\]
which is defined by
\[
HD(f)\Lambda =D_{f\mu,\mu'(f \times f)}
\]
where $\Lambda:X\times X \to X\wedge X$ is the natural quotient
map.
\end{defn}
\begin{rem}
The definition of H-derivation depends on the H-space structures
on both $X$ and $Y$.
\end{rem}
%\begin{rem}
%It is easy to know that H-derivation give a homomorphism of groups
%if the multiplication on the target space is homotopy associative
%and homotopy abelian .
%\end{rem}
\begin{rem}
Let $(X,\mu)$ and  $(Y,\mu')$ be H-spaces and $f: X \to Y$ be a
map of spaces. It is well known that $f:(X,\mu_X) \to (Y,\mu_Y)$
is an H-map
 iff $HD(f)=*$.
\end{rem}

An easy but crucial corollary of this last remark is  the
following
\begin{cor}\label{T:h-map}
Let  $(X,\mu)$ and  $(X',\mu')$ be H-spaces and $f: X \to X'$ be a
map of spaces. Assume further that $\pi_iX'=0$ for $i\geq 2d$ and
$X$ is $(d-1)$-connected . Then $f$ is a $\mu -\mu'$ H-map.
\end{cor}
Following is one of the fundamental properties of H-derivation
\begin{prop}\label{T:deri}
Let $(X_i,\mu_i)$ be H-spaces, $i=0,1,2$ and $f:X_0 \to X_1$,
$g:X_1 \to X_2$ be maps of spaces. \newline (a)If $f:X_0 \to X_1$
is a $\mu_0 -\mu_1$ H-map, then $HD(gf)=HD(g)(f\wedge f)$ \newline
(b)If $g:X_1 \to X_2$ is a $\mu_1 -\mu_2$ H-map, then
$HD(gf)=gHD(f)$
\end{prop}

Another ingredient for the main result is the phantom map . Recall
that a map $f$ from a CW complex $X$ is called an phantom map if
its restriction to the $n$-th skeleton is inessential for any
integer $n$.  Let $Ph(X,Y)$ denote the set of homotopy classes of
phantom maps from $X$ to $Y$. The following result which follows
from the Sullivan conjecture provides us many examples of phantom
maps.
\begin{thm}\cite{cm95}\label{T:phantom}
Let $Y=\Omega^i K$ and $X=\Sigma^j Z$ such that $i,j\geq 0$ ,$K$
is a 1-connected finite CW complex.  Then every map from $X$ to
$Y$ is a phantom map if $Z$ is as follows:
\begin{itemize}
\item{$Z$ is the classifying space of a 1-connected compact Lie group}
\item{$Z$ is an infinite loop space with torsion fundamental group}
\item{$Z$ has only finitely number of nontrivial homotopy groups}
\end{itemize}
and in this case we have
\[
Ph(X,Y)=[X,Y]=[X_{(0)},Y]=\prod_{n>0}H^n(X,\pi_{n+1}(Y)\bigotimes
R)
\]
\end{thm}
Let $Ph(X,Y)$ denote the set of homotopy classes of phantom maps
from $X$ to $Y$.  The p-localization $l_p$ induces a natural map
\[
l_p^*:Ph(X,Y) \to Ph(X,Y_{(p)})
\]
It follows that there is a natural map
\[
l:Ph(X,Y) \to \prod_p Ph(X,Y_{(p)})
\]
It is well known that $l$ is an epimorphism \cite{steiner} and
$Ker(l)$ is nontrivial iff $Ph(X,Y)$ is nontrivial\cite{mcgib},
see also \cite{roitharp},\cite{roit2} . The phantom map in
$Ker(l)$ is called special , following Roitberg\cite{roit}, see
also,\cite{moller} where it is called the clone of constant map.

On the other hand , Gray, Le Minh Ha , McGibbon and Strom
\cite{gray}, \cite{lestrom},\cite{mstrom} introduced the notion of
Gray index which is defined as follows:
\begin{defn}
Let $f:X \to Y$ be a phantom map. Then  $f$ can be factorized as
the composition $X \to X/X_k \overset{\bar{f}}{\to} Y$ for each
$k$ . The Gray index of $f$ , denoted by $G(f)$, is the largest
integer $k$ such that the $\bar{f}$ can be chosen to be a phantom
map. $G(f)=\infty$ if no such $k$ exists.
\end{defn}
\begin{rem}\label{T:dindex}
Let $f:X \to Y$ be a phantom map. Then  $f$ can be lifted to the
k-th connected covering for each $k$ and $G(f)+1$ is the largest
integer $k$ such that the $k$th lifting can be chosen to be a
phantom map.
\end{rem}
A useful fact we need is
\begin{prop}\label{T:index}
Let  $f:X \to Y$ be a phantom map. Then \newline (i)$G(f)\geq n$
if $X$ is n-connected or $Y$ is $n+1$-connected. \newline (ii)
 $G(f)\in \{k | H^n(X,\pi_{n+1}(Y)\bigotimes
Q)\neq 0 \}$ if $H^n(X,\pi_{n+1}(Y)\bigotimes Q)=0$ for n
sufficiently large.
\end{prop}
For the proof of the Proposition above, see \cite{lestrom} and
\cite{mstrom}.

An immediate corollary of the above Proposition which is crucial
to our purpose is
\begin{cor}\label{T:cindex}
Let $f:K(Z,2m)\wedge K(Z,2m) \to \Omega^1 S^{4m+2}$ be any
essential  map. Then $G(f)=4m$ where $m \geq 1$.
\end{cor}
 Now we are ready to  prove the main result.
\section{Proof of Theorem\ref{T:main}}
First is a preliminary lemma needed later.
\begin{lem}\label{T:prelim}
Let $X$ be an H-space which is $(e-1)$-connected and $\pi_i(X)=0$
for $i\geq 2e$ and $Y=\Omega^j K$ where $K$ is a
$(d+j-1)$-connected finite CW-complex with $j\geq 1$ and $d> e\geq
2 $. Let $\psi: X \to Y$ be any essential map. Then $Z$(=homotopy
fiber of $\psi$) is not an H-space if $HD(\psi)\circ (i \wedge i)$
is essential where $i:Z \to X$ is the homotopy fiber of $\psi$ .
\end{lem}

\begin{proof}
If $Z$ is an H-space , then $*= \psi \circ i:Z \to Y $ is an H-map
and  $HD(*)=0$. Since $i$ is an H-map by Corollary\ref{T:h-map} ,
it follows by Proposition\ref{T:deri} that $$*=HD(*)=HD(\psi \circ
i)=HD(\psi)\circ (i \wedge i)$$ which is in contradiction to the
condition.
\end{proof}
\begin{rem}
To apply the above lemma it suffices to discuss  when $\psi$ is
not an H-map  and when $i\wedge i$ induces an injective.
\end{rem}
\begin{thm}\label{T:loop}
Let $X=K(\mathbb{Z},2m)$ and $Y=\Omega^1 S^{4m+2}$ with $m\geq 1$.
Then there is no essential H-map from $X$ to $Y$.
\end{thm}
\begin{proof}
By Proposition\ref{T:phantom}, the rationalization $r:X \to
X_{(0)}$ which is an H-map induces an isomorphism of groups
\[
r^*:[X_{(0)},Y]\to [X,Y]
\]
It follows from Proposition\ref{T:deri} that it suffices to prove
that  there is no essential H-map from $X_{(0)}$ to $Y$.

On the other hand the map $h:S^{4m+2}\to K(Z,4m+2)$ which
represents a generator of $H^{4m+2}(S^{4m+2};\mathbb{Z})\cong
\mathbb{Z}$ induces an isomorphism of groups
\[
(\Omega^3 h)_*:[X_{(0)},Y] \to [X_{(0)},K(\mathbb{Z},4m+1)]
\]
Again the Proposition\ref{T:deri} implies  that it suffices to
prove that  there is no essential H-map from $X_{(0)}$ to
$K(\mathbb{Z},4m+1)$ which is well known to be equivalent to the
injectivity of the following homomorphism
\[
\theta^*:H^{4m+2}(\Sigma X_{(0)};\mathbb{Z})\to H^{4m+2}(\Sigma
X_{(0)}\wedge X_{(0)};\mathbb{Z})
\]
where $\theta$ is defined as follows:

Let $X*X= X\times I\times X /\{(x,0,y)\sim (x,0,y'),(x,1,y)\sim
(x',1,y)\}$ be the join . There is a well defined map
\[
k:X*X \to \Sigma X\wedge X
\]
by k[x,t,y]=(x,y,t). It is well known that $X*X$ is homotopy
equivalent to $\Sigma X\wedge X$. If $X$ is an H-space with
multiplication $\mu$ , then $\theta$ is the composite map
\[
\Sigma X\wedge X \simeq X*X \overset{k}{\to} \Sigma (X\times X)
\to \Sigma X
\]
where the last map is the map $-\Sigma \pi_1 + \Sigma \mu -\Sigma
\pi_2 $ and $\pi_1$, $\pi_2$ are the projection of $X\times X$ to
the first and second factors respectively .

Consider the following commutative diagram where the horizontal
maps  which are isomorphisms come from the universal coefficient
Theorem
\[
\begin{CD}
H^{4m+2}(\Sigma X_{(0)};\mathbb{Z}) @>>>
\text{Ext}(H_{4m+1}(\Sigma X_{(0)};\mathbb{Z}),\mathbb{Z}) \\
 @V{\theta}^*VV @V\text{Ext}({\theta}_*,\mathbb{Z})VV \\
H^{4m+2}(\Sigma X_{(0)}\wedge X_{(0)};\mathbb{Z}) @>>>
\text{Ext}(H_{4m+1}(\Sigma X_{(0)}\wedge
X_{(0)};\mathbb{Z}),\mathbb{Z})
\end{CD}
\]
It follows that it suffices to prove that the map
\[
\theta_*:H_{4m+1}(\Sigma X_{(0)}\wedge X_{(0)};\mathbb{Z}) \to
H_{4m+1}(\Sigma X_{(0)};\mathbb{Z})
\]
or equivalently the map
\[
\theta_*:H_{4m}(X\wedge X;\mathbb{Q}) \to H_{4m}( X;\mathbb{Q})
\]
is  injective . On the other hand, it is well known that the map
$\theta_*$ is dual to the reduced coproduct which is an
isomorphism in this case and thus completes the proof.
\end{proof}

The Theorem\ref{T:main} is actually a corollary of the following
more general Theorem
\begin{thm}
Let $X=K(\mathbb{Z},2m)$ and $Y=\Omega S^{4m+2}$ with $m\geq 1$.
Let $\psi: X \to Y$ be any essential, special, phantom map. Then
$Z$(=homotopy fiber of $\psi$) is not an H-space and has the same
Mislin genus with $X \times \Omega Y$.
\end{thm}
\begin{proof}
That $Z$ and $X \times \Omega Y$ have the same Mislin genus
follows from the condition that $\psi: X \to Y$ is a  special
phantom map.

On the other hand Lemma\ref{T:prelim} and Theorem\ref{T:loop}
apply here. Thus the Theorem above follows from the following
Proposition.
\end{proof}
\begin{prop}
Let $X=K(\mathbb{Z},2m)$ and $Y=\Omega S^{4m+2}$ with $m\geq 1$.
Let $\psi: X \to Y$ be any essential, special, phantom map which
exists by Proposition\ref{T:phantom} and the remark after it. Then
$(i \wedge i)^*:[X\wedge X,Y]\to [Z\wedge Z,Y]$ is injective where
$i:Z \to X$ is the homotopy fiber of $\psi$ .
\end{prop}
\begin{proof}
Let $f:X \wedge X\to Y$ be any essential map . If $f \circ (i
\wedge i) \sim *$ we will prove that this leads to a contradiction
which concludes the proof. Since $f$ is a phantom map , $f \circ
(i \wedge i) \sim *$ is also a phantom map.Thus we have the
following commutative diagram up to homotopy.
\[
\begin{CD}
Z \wedge Z @>i\wedge i>> X\wedge X @>f>> Y \\ @VVV @VVV @VidVV \\
Z_{(0)}\wedge Z_{(0)} @>i_{(0)}\wedge i_{(0)}>>  X_{(0)}\wedge
X_{(0)} @>\tilde{f}>> Y
\end{CD}
\]
 If $f \circ (i
\wedge i) \sim *$ , then $\tilde{f}\circ (i_{(0)}\wedge i_{(0)})$
is the composite  $Z_{(0)}\wedge Z_{(0)} \to \Sigma Z_{\tau}\wedge
 Z_{\tau} \overset{h}{\to} Y$ where $Z_{\tau}$ is the
homotopy fiber of the rationalization $X \to X_{(0)}$. On the
other hand we claim that
\begin{claim}\label{T:claim}
Any map $h:\Sigma Z_{\tau}\wedge  Z_{\tau} \to Y$ factors through
a map $\Sigma Z_{\tau}\wedge  Z_{\tau} \to \Sigma F_{\tau}\wedge
 F_{\tau} $ where $F=\Omega^2 S^{4m+2}$.
\end{claim}
Assuming this , note that $i_{(0)}\wedge i_{(0)}$ admits a right
inverse, we have that $f$ is the composite $X \wedge X \to
X_{(0)}\wedge X_{(0)} \to \Sigma F_{\tau}\wedge  F_{\tau} \to Y$.

It is easy to know that $\Sigma F_{\tau}$ is $4m-1$-connected and
thus $ \Sigma F_{\tau}\wedge  F_{\tau}$ is $8m-2$-connected. It
follows that $f$ is the composite $X \wedge X \to X_{(0)}\wedge
X_{(0)} \to Y<8m-2> \to Y$. By Remark\ref{T:dindex}, $G(f)\geq
8m-3$ which contradicts Corollary\ref{T:cindex}.
\end{proof}
\begin{rem}
Roitberg has shown us how the use of Gray index can be avoided.
\end{rem}
It remains to prove the Claim\ref{T:claim} which follows from the
following
\begin{lem}
There is a map $\Sigma g\wedge g:\Sigma Z_{\tau}\wedge  Z_{\tau}
\to \Sigma F_{\tau}\wedge F_{\tau}$ such that the following map is
a weak homotopy equivalence
\[
(\Sigma g\wedge g)^*:map_*(\Sigma F_{\tau}\wedge  F_{\tau},Y)\to
map_*(\Sigma Z_{\tau}\wedge \ Z_{\tau},Y)
\]
\end{lem}
\begin{proof}
Roitberg and P. Touhey  proved in \cite{touhey} that , if $X,Y$
have the same Mislin genus, then $X_{\tau}\simeq Y_{\tau}$. So we
have a map $g$ which is a composite $$ Z_{\tau} \simeq (X\times
F)_{\tau} \overset{ \pi_{\tau}}{\to} F_{\tau}$$ where $\pi:
X\times F \to F$ is the projection.

To prove $\Sigma g \wedge g$ induces a homotopy equivalence it
suffices to prove that $ \pi_{\tau}: (X\times F)_{\tau} \to
(F)_{\tau}$ induces a homotopy equivalence
\[
(\Sigma \pi_{\tau}\wedge \pi_{\tau})^*:map_*(\Sigma F_{\tau}\wedge
F_{\tau},Y)\to map_*(\Sigma (X\times F)_{\tau}\wedge (X\times
F)_{\tau},Y)
\]
which follows directly from the fact that $map_*(X_{\tau},Y)$ is
weakly contractible and the fact  $$map_*(\Sigma X_{\tau},Y)
\simeq map_*(\Sigma X,\hat{Y})$$ which can be found in
\cite{touhey1} , for a stronger result , see Pan and Woo
\cite{panwoo}.
\end{proof}
----------------------------------------------------------------

----------------------------------------------------------------

\end{document}